%%%%%%     AMS-TeX     %%%%%%

\documentstyle{amsppt}
\hyphenation{semi-ring non-emp-ty}
\TagsOnRight
\NoRunningHeads

%%%%%%%%%%%%%%%%%%%%%%%%%%%%%%%%%%%%%%%%%%%%%

\def \!-#1{$#1\text -$\nolinebreak}
\def\!inf-{$\wedge \text -$\nolinebreak}

\def \inn{\nomathbreak\in\nomathbreak}
\def \ssubset{\nomathbreak\subset\nomathbreak}

%\define \Sup{\operatornamewithlimits {sup} \limits}

\define \OSup{\operatornamewithlimits \bigoplus \limits}

\def\<{\preccurlyeq }

\def\>{\succcurlyeq }

\define \0{{\bold 0}}
\define \1{{\bold 1}}

\def\a-{$a\text -$\nolinebreak }
\def\aa-{$a\text -$\nolinebreak }
\def\b-{$b$-\hskip0pt}
\def\c-{$\vee\text -$\nolinebreak }
\def\w-{$\wedge$-\hskip0pt}
\def\d-{$\delta$-\hskip0pt}
\def\oo-{$o\text -$\nolinebreak }
\def\wo-{$wo\text -$\nolinebreak }
\def\G-{$G\text -$\nolinebreak }
\def\g-{$g\text -$\nolinebreak }

\define \R-{{${\bold R_{\oplus }}\text -$\nolinebreak }}
\define \plus-{{$\oplus\text -$\nolinebreak}}

\define \Rmax{{\bold R}_{\max}}

\define \disk#1#2{D_{#1} (#2)}
\define\bR{\text{\bf R}}

\def\krn{\text{\it kern\/}}
\def\id{\text{\it id\/}}
\def\Id{\text{\it Id\/}}

%%%%%%%%%%%%%%%%%%%%%%%%%%%%%%%%%%%%%%%%%%%%%
\topmatter
\title  Nuclear semimodules and kernel theorems
in idempotent analysis. An algebraic approach
\endtitle
\author G.~L.~Litvinov and G.~B.~Shpiz\endauthor
\endtopmatter
%%%%%%%%%%%%%%%%%%%%%%%%%%%%%%%%%%%%%%%%%%%%%

\document
In this note we describe conditions under which, in idempotent
functional analysis (see \cite{1--3}), linear operators have integral
representations in terms of idempotent integral of V.~P.~Maslov. We
define the notion of nuclear idempotent semimodule and describe
idempotent analogs of the classical kernel theorems of L.~Schwartz and
A.~Grothendieck (see, e.g., \cite{4--6}).  In \cite{11}, for the
idempotent semimodule of all bounded functions with values in the
Max-Plus algebra, there was posed a problem of describing the class of
subsemimodules where some kind of kernel theorem holds.  Some of the
results obtained here can be regarded as possible versions of an
answer to this question. Previously, some theorems on integral
representations were obtained for a number of specific semimodules
consisting of continuous or bounded functions taking values mostly in
the Max-Plus algebra (see, e.g., \cite{7--10, 1, 3}).

In this work (and further publications), a very general case of
semimodules over boundedly complete idempotent semirings is considered.
This note continues the series of publications \cite{1--3}; we use the
notation and terminology defined in those articles.

{\bf 1. Functional semimodules.} Idempotent analysis is based on replacing
the number fields by idempotent semifields and semirings. In other words, a
new set of basic associative operations (a new addition~$\oplus$ and new
multiplication~$\odot$) replaces the traditional arithmetic
operations. These new operations must satisfy all axioms of semifield or
semiring; also, the new addition must be idempotent, i.e., satisfy $x
\oplus x = x$ for all~$x$ belonging to the semifield or semiring (see,
e.g., \cite{7--11}). A typical example is the semifield
$\Rmax=\bR\bigcup\{-\infty\}$, which is known as the Max-Plus algebra. This
semifield consists of all real numbers and an additional element~$-\infty$,
denoted by~$\0$, that is the zero element of $\Rmax$; the operations are
defined by $x\oplus y = \max\{x,y\}$ and $x\odot y = x+y$. The unit
element~$\1$ of~$\Rmax$ coincides with the usual zero. One can find a
lot of nontrivial examples of idempotent semirings and semifields, e.g., in
\cite{3, 7--11}.  An {\it idempotent semimodule\/} over an idempotent
semiring~$K$ is an additive commutative idempotent semigroup, with the
addition operation denoted by~$\oplus$ and the zero element denoted
by~$\0$, such that a product $k \odot x$ is defined, for all~$k$ in~$K$
and $x$ in the semimodule, in such a way that the usual rules are
satisfied.  For exact definitions (including those of $b$-linearity and
$b$-completeness) and examples, see \cite{1--3}. Recall that any idempotent
semigroup has the canonical partial order: $x\<y$ iff $x\oplus y = y$.

A~subset of an idempotent semimodule is called a~{\it subsemimodule\/} if
it is closed under addition and multiplication by scalar coefficients.
A~subsemimodule of a~\b-complete semimodule is called a {\it
\b-subsemimodule\/} if the corresponding embedding map is a
\b-homomorphism. The main feature of \b-subsemimodules is that
restrictions of \b-linear operators and functionals to these semimodules
are \b-linear.

Let $X$ be an arbitrary nonempty set and $K$ an idempotent semiring.
By~$K(X)$ denote the semimodule of all maps $X \to K$ endowed with
pointwise operations. If $K$ is a \b-complete semiring, then $K(X)$ is a
\b-complete semimodule. For a point~$x \inn X$, by $\delta_x$ denote the
functional on $K(X)$ taking $f$ to $f(x)$. It can be easily checked that
the functional $\delta_x$ is \b-linear on $K(X)$.

We call an {\it idempotent functional semimodule} on a set~$X$, or an
IFS for short, a subset of $K(X)$ that is invariant under
multiplication by scalars and is an upper semilattice with respect to
the induced order. An IFS will always be regarded as a semimodule with
respect to the induced multiplication by coefficients from~$K$ and the
addition defined by $x\oplus y = \sup\{x, y\}$, where $\sup$ is taken
inside the functional semimodule. An IFS is called {\it \b-complete\/}
if it is a \b-complete semimodule. An IFS on the set~$X$ is called a
{\it functional \w-semimodule\/} if it contains $\0$ and is closed
under the operation of taking infima of nonempty subsets in
$K(X)$. Evidently, any functional \w-semimodule is a \b-complete
semimodule. Notice that a functional semimodule (and even a
\w-semimodule) may not be a subsemimodule of $K(X)$, because it
inherits from~$K(X)$ the order but, in general, not the addition
operation.

In general, a functional of the form $\delta_x$ on a functional semimodule
is not even linear, much less \b-linear. On the other hand, the following
proposition holds, which is a direct consequence of the definitions.

\proclaim{Proposition 1} A \b-complete functional semimodule~$W$
on a set~$X$ is a \b-subsemimodule of $K(X)$ iff each functional of the form
$\delta_x$ (with $x \inn X$) is \b-linear on~$W$.\endproclaim

{\bf 2. Integral representations of linear operators in functional
semimodules.}
In what follows, $K$ always denotes a \b-complete semiring unless otherwise
specified. Let $W$ be an idempotent \b-complete semimodule over~$K$ and $V$
a subset of $K(X)$ that is a \b-complete IFS on~$X$. A map $A\colon V
\to W$ is called an {\it integral\/} operator or an operator with an
{\it integral representation\/} if there exists a map $k \colon X \to W$,
called the {\it kernel of the operator~$A$\/}, such that
$$
Af=\OSup_{x\in X} f(x)\odot k(x). \tag 1
$$
In idempotent analysis, the right-hand side of formula~\thetag{1} is often
written as $\int_X^{\oplus}f(x)\odot k(x)\, dx$ (see \cite{7--9}).
Regarding the kernel~$k$ it is supposed that the set $\{\,f(x)\odot k(x)
\mid x\inn X\,\}$ is bounded in~$W$ for all $f \inn V$ and $x \inn X$. We
denote the set of all functions with this property by $\krn_{V,W}(X)$. In
particular, if $W = K$ and $A$ is a functional, then this
functional is called {\it integral}.

We call a functional semimodule $V \subset K(X)$ {\it
non-degenerate\/} if for each point~$x \inn X$ there exists a
function~$g \inn V$ such that $g(x) = \1$, and {\it admissible\/} if
for each function~$f \inn V$ and each point~$x \inn X$ there exists a
function~$g \inn V$ such that $g(x) = \1$ and $f(x) \odot g \<
f$. Note that all idempotent functional semimodules over a semifield
are admissible (it suffices to set $g = f(x)^{-1} \odot f$).

If an operator has an integral representation, this representation may not
be unique.  However, if the semimodule~$V$ is non-degenerate, then the set
of all kernels of a fixed integral operator is bounded with respect to the
natural order in the set of all kernels and is closed under the supremum
operation applied to its arbitrary subsets. In particular, any integral
operator defined on a non-degenerate IFS has a unique maximal kernel.
An important point is that an integral operator is not necessarily
\b-linear and even linear except when $V$ is a \b-subsemimodule of~$K(X)$
(see Proposition~2 below).

If $W$ is an IFS on a nonempty set $Y$, then the kernel~$k$ of an
operator~$A$ can be naturally identified with a function on $X \times Y$
defined by $k(x,y) = (k(x))(y)$. This function will also be called the
{\it kernel\/} of the operator~$A$. As a result, the set $\krn_{V,W}(X)$ is
identified with the set $\krn_{V,W}(X,Y)$ of all maps $k\colon X \times Y
\to K$ such that for each~$x \inn X$ the map $k_x\colon y \mapsto k(x, y)$
lies in~$W$ and for each~$v \inn V$ the set $\{\, v(x)\odot k_x \mid x \inn
X \,\}$ is bounded in~$W$. Accordingly, the set of all kernels of \b-linear
operators can be embedded in $\krn_{V,W}(X,Y)$.

If $V$ and~$W$ are \b-subsemimodules of~$K(X)$ and of~$K(Y)$, respectively,
then the set of all kernels of \b-linear operators from~$V$ to~$W$ can be
identified with $\krn_{V,W}(X,Y)$ (see Proposition~2 below), and the
following formula holds:
$$
Af(y)=\OSup_{x\in X} f(x)\odot k(x,y)=\int_X^{\oplus} f(x)
\odot k(x,y) dx. \tag 2
$$
This formula coincides with the usual definition of an operator's integral
representation. Note that formula~\thetag{1} can be rewritten in the form
$$
Af=\OSup_{x\in X} \delta _x(f)\odot k(x). \tag 3
$$

\proclaim{Proposition 2} A \b-complete functional semimodule~$V$ on a
nonempty set~$X$ is a \b-subsemimodule of~$K(X)$ iff all integral operators
defined on~$V$ are \b-linear.\endproclaim

\proclaim{Theorem 1} Suppose an admissible \b-complete semimodule~$W$ over
a \b-complete semiring~$K$ and a functional \w-semimodule~$V \ssubset K(X)$
are fixed. Then each \b-linear operator $A\colon V \to W$ has
an integral representation of the form~\thetag{1}. In particular, if $W$ is
a \b-subsemimodule of $K(Y)$, then the operator~$A$ has an integral
representation of the form~\thetag{2}.\endproclaim

{\bf 3. Integral representations of \b-nuclear operators.}
We now introduce several important definitions. A map between \b-complete
semimodules $g\colon V \to W$ is called {\it one-dimensional\/} (or a {\it
rank\/}~1 map) if it has the form $v \mapsto \phi(v) \odot w$, where $\phi$
is a \b-linear functional on~$V$ and $w$ lies in~$W$. A~map~$g$ is called
{\it \b-nuclear\/} if it is a sum of a bounded set of one-dimensional
maps. It goes without saying that \b-nuclear maps are closely related to
tensor products of idempotent semimodules, see~\cite{2}.

\proclaim{Theorem 2} Let $W$ be a \b-complete semimodule over a \b-complete
semiring~$K$. If all \b-linear functionals defined on a \b-subsemimodule~$V$
of~$K(X)$ are integral, then a \b-linear operator $A\colon V \to W$ has an
integral representation iff it is \b-nuclear.\endproclaim

{\bf 4. The \b-approximation property and \b-nuclear semimodules.} We
shall say that a \b-complete semimodule~$V$ has the {\it
\b-approximation property\/} if the identity operator $\id\colon V \to
V$ is \b-nuclear (for a~treatment of the approximation property of
locally convex spaces in the traditional analysis, see \cite{4, 5}).

Let $V$ be an arbitrary \b-complete semimodule over a \b-complete
idempotent semiring~$K$. We call this semimodule a {\it \b-nuclear
semimodule\/} if any \b-linear map of~$V$ to an arbitrary \b-complete
semimodule~$W$ over~$K$ is a \b-nuclear operator. Recall that, in the
traditional analysis, a locally convex space is nuclear iff all continuous
linear maps of this space to any Banach space are nuclear operators,
see~\cite{4, 5}.

\proclaim{Proposition~3} Let $V$ be an arbirary \b-complete semimodule over
a \b-complete semiring~$K$. The following statements are equivalent:
\roster
\item the semimodule~$V$ has the \b-approximation property;
\item each \b-linear map from~$V$ to an arbitrary \b-complete
semimodule~$W$ over~$K$ is \b-nuclear;
\item each \b-linear map from an arbitrary \b-complete semimodule $W$
over~$K$ to the semimodule~$V$ is \b-nuclear.
\endroster
\endproclaim
\proclaim{Corollary} A \b-complete semimodule over a \b-complete
semiring $K$ is \b-nuclear iff this semimodule has the \b-approximation
property.\endproclaim
Recall that, in the traditional analysis, any nuclear space has the
approximation property but the converse statement is not true.

{\bf 5. Kernel theorems in function semimodules.}
Let a semimodule~$V$ of~$K(X)$ be a \b-complete IFS over a
\b-complete semiring~$K$. We shall say that the {\it kernel theorem\/} holds
in the semimodule~$V$ if all \b-linear maps of this semimodule to an
arbitrary \b-complete semimodule over~$K$ have an integral representation.

\proclaim{Theorem 3} Suppose a \b-complete semiring~$K$ and a nonempty
set~$X$ are given. The kernel theorem holds in a \b-subsemimodule~$V$
of~$K(X)$ iff all \b-linear functionals on~$V$ have integral
representations and the semimodule~$V$ is \b-nuclear, i.e., has the
\b-approximation property.\endproclaim

\proclaim{Theorem 3a} Suppose a \b-complete semiring~$K$ and a nonempty
set~$X$ are given. The kernel theorem holds in a \b-subsemimodule~$V$
of~$K(X)$ iff the identity operator $\id\colon V \to V$ is
integral.\endproclaim

{\bf 6. Integral representations of operators in abstract idempotent
semimodules.}
Let $V$ be a \b-complete idempotent semimodule over a \b-complete
semiring~$K$ and $\phi$~be a \b-linear functional defined on~$V$. We call
this functional a {\it \d-functional\/} if there exists an element $v \in
V$ such that $\phi(w) \odot v \< w$ for all $w\inn V$. Denote
by~$\Delta(V)$ the set of all \d-functionals on~$V$ and
by~$i_\Delta$ the natural map $V \to K(\Delta(V))$ given by $(i_\Delta
(v))(\phi) = \phi(v)$ for all $\phi \inn \Delta(V)$. Note that if a
one-dimensional operator~$\phi \odot v$ appears in a decomposition of the
identity operator on~$V$ into a sum of one-dimensional operators, then
$\phi \inn \Delta(V)$.

Denote by $\id$ and~$\Id$ the identity maps on the semimodules $V$
and~$i_\Delta(V)$, respectively.

\proclaim{Proposition 4} The operator~$\id$ is \b-nuclear iff $i_\Delta$ is
an embedding and the operator~$\Id$ is integral.\endproclaim

\proclaim{Theorem 4} A \b-complete idempotent semimodule~$V$ over a
\b-complete semiring~$K$ is isomorphic to an IFS in which the kernel
theorem holds iff the identity map on~$V$ is a \b-nuclear operator, i.e.,
iff $V$ is a \b-nuclear semimodule.\endproclaim

The following statement shows that, in a certain sense, the
embedding~$i_\Delta$ is a universal representation of a \b-nuclear
semimodule in the form of an IFS in which the kernel theorem holds.

\proclaim{Proposition 5} Let $K$ be a \b-complete semiring, $X$ a nonempty
set and $V \ssubset K(X)$ a functional \b-semimodule on~$X$ in which the
kernel theorem holds. Then there exists a natural map $i\colon X \to
\Delta(V)$ such that the corresponding map $i_*\colon K(\Delta(V)) \to
K(X)$ maps $i_\Delta(V)$ on~$V$ isomorphically.\endproclaim

The authors are cincerely grateful to V. N. Kolokoltsov, V. P. Maslov,
and A. N. Sobolevski\u\i\ for valuable suggestions and support.
The work was supported by RFBR under Grant~99--01--00196 and by the Dutch
Organization for Scientific Research N.W.O.

%\newpage


\Refs

\ref\no 1
\by {\it Litvinov G. L., Maslov V. P., Shpiz G. B.}
Linear functionals on idempotent spaces:
An algebraic approach // Doklady Akademii Nauk.~--- 1998.~--- V.~363.,
N.~3.~--- pp.~298--300 (in Russian). See also an English translation in:
Doklady Math.~--- 1998.~--- V.~58, N.~3.~--- pp.~389--391 \endref

\ref\no 2
\by {\it Litvinov G. L., Maslov V. P., Shpiz G. B.}
Tensor products of idempotent semimodules.
An algebraic approach //
Matematicheskie Zametki.~--- 1999.~--- V.~65, N.~4.~--- pp.~479--489 (in
Russian). See also an English translation in: Math. Notes.~--- 1999.~---
V.~65, N.~4.~--- pp.~572--585
\endref

\ref\no 3
\by {\it Litvinov G. L., Maslov V. P., Shpiz G. B.}
Idempotent functional analysis: An algebraic approach //
Matematicheskie Zametki.~--- 2001.~--- V.~69, N.~5.~--- pp.~696--729 (in
Russian). See also an English translation in: Math. Notes.~-- 2001.~---
V.~69, N.~5.~--- pp.~758--797
\endref

\ref\no 4
\by {\it Grothendieck A.}
Produits tensoriels topologiques et espaces nucl\'eaires.--
Mem. Amer. Math. Soc., V.~16, Providence, R.I., 1955
\endref

\ref\no 5
\by {\it Schaefer H. H.}
Topological vector spaces.~---
New York: The Macmillan Company and London: Collier-Macmillan Ltd, 1966
\endref

\ref\no 6
\by {\it Schwartz L.}
Th\'eorie des noyaux.~---
Proc. of the Intern. Congress of Math.~--- 1950.~--- V.~1.~--- pp.~220--230
\endref

\ref\no 7
\by  {\it Maslov V.P.}
M\'ethodes op\'eratorielles.~--- Moscow: Mir, 1987
\endref

\ref\no 8
\by {\it Maslov V. P., Kolokoltsov V. N.}
Idempotent analysis and its application in optimal control.~---
Moscow: Nauka, 1994 (in Russian)
\endref

\ref\no 9
\by {\it Maslov V.P., Samborski\u\i\ S.N., Eds.}
Idempotent Analysis.~--- Adv. in Sov. Math., V.~13.~--- Providence, R.I.:
Amer. Math. Soc., 1992 \endref

\ref\no 10
\by {\it Shubin M. A.}
Algebraic remarks on idempotent semirings
and a kernel theorem in spaces of bounded functions.~---
Moscow: The Institute for New Technologies, 1990 (in Russian)
\endref

\ref\no 11
\by {\it Gunawardena J.}
An introduction to idempotency.~--- Idempotency/J. Guna\-wa\-rdena, Ed.~---
Cambridge: Cambridge University Press, 1998.~--- pp.~1--49 \endref

\endRefs
\enddocument